\documentclass[11pt]{article}
\usepackage{graphics,amsmath,amssymb,float}
\usepackage{latexsym}
\usepackage{epsfig}
\usepackage{hyperref}
\providecommand{\abs}[1]{\left\lvert#1\right\rvert}

\newcommand{\set}[1]{\left\{#1\right\}}
\def\epsilon{\varepsilon}

\def\phi{\varphi}

\def\om{\omega}

\newtheorem{theorem}{Theorem}[section]
\newtheorem{lemma}[theorem]{Lemma}

\newtheorem{proposition}[theorem]{Proposition}
\newtheorem{remark}[theorem]{Remark}

\def\N{{\mathbb N}}
\def\C{{\mathbb C}}
\def\R{{\mathbb R}}

\newenvironment{Proof}{\removelastskip\par\medskip
\noindent{\em Proof.} \rm}{\penalty-20\null$\square$\par\medbreak}

\title{\bf  Inverse Observability Inequalities \\ for Integrodifferential Equations \\ in Square Domains }

\author{Paola Loreti
\thanks{Dipartimento di Scienze di Base e Applicate per l'Ingegneria
Sapienza Universit\`a di Roma,
Via Antonio Scarpa 16, 00161 Roma (Italy); e-mail: 
$<$paola.loreti@uniroma1.it$>$ }
\and Daniela Sforza
\thanks{Dipartimento di Scienze di Base e Applicate per l'Ingegneria,
 Sapienza Universit\`a di Roma,
Via Antonio Scarpa 16, 00161 Roma (Italy); e-mail: 
$<$daniela.sforza@uniroma1.it$>$ }}

\begin{document}
\date{}

\maketitle

\begin{abstract}
In this paper  we will consider oscillations of square viscoelastic membranes by
adding to the wave equation another term, which takes into account the memory.
To this end, we will  study a class of integrodifferential equations in  square domains.
By using  accurate estimates of the spectral properties of the integrodifferential operator,
we will prove an inverse observability inequality.
\end{abstract}

\bigskip
\noindent
{\bf Keywords:} observability; Fourier series; 
Ingham estimates

\smallskip
\noindent
{\bf MSC:   45K05} 
\section{Introduction}
In \cite {LS0} and \cite {LS} we solved  a Dirichlet boundary control problem
for the  wave equation with  the exponential  memory  kernel
$$k(t)=\beta e^{-\eta t}. $$
 The result was established under some conditions on the parameters $\beta$ and
$\eta$, that is  $\eta>3\beta/2$, $\eta \geq 0$, \,\, $\beta\geq 0$. 
In the investigation a key point to get an estimate for the control time was to prove the inverse inequality.
\noindent
In \cite {LS} the analysis was done in the one-dimensional case, obtaining a precise estimate of the observability time $T$.  In \cite{LoretiSforza2}  we also solved the problem for $n$-dimensional balls.
It was an open problem to extend to simple  domains like rectangles, common  in applications, the previous results using the Fourier method. 
The inverse observability estimate for the wave equation without memory  was obtained
under the geometrical  condition that the control time $T$ is greater than twice the diagonal of the rectangle,
see \cite{BLR}.
By means of the Fourier method, Mehernberger \cite {Meh2009} obtained  a weaker result,  nevertheless his method
can be adapted  to get the inverse inequality for other models. See also   \cite {KomLor159}  for an improvement of \cite {Meh2009} and further applications.

For another approach we have to mention the paper \cite{K2}, where the analysis of the kernel is done by a compact perturbation and 
the author proves his result by means of
a unique continuation property of the integro-differential equation. The method proposed by us is more direct and may be easily adapted to other boundary conditions, since in  our estimates  the dependence
on  the eigenvalues of the integro-differential operator is explicitly given.
Moreover, Theorem 1.3 below has an interest in itself, because it contains a method
 which is more general than those used in  \cite{Meh2009,KomLor159}.

It is noteworthy that exponential kernels arise in  viscoelasticity theory, such as in the analysis of Maxwell fluids  or Poynting -Thomson solids,
see  e.g.
\cite {Pruss,Re1}.  For other references in viscoelasticity theory see the seminal papers of Dafermos \cite{D1,D2} and  \cite{RHN,LPC}.

In this paper  we will consider oscillations of square viscoelastic membranes by
adding to the wave equation another term, which takes into account the memory.
We will fix $\eta=3\beta/2$ 
to study the integrodifferential equation in a square.
This assumption has the double target to  simplify the computation for the square and to extend to the $2$-d case the results given in  
\cite{LS}. 

We will go back to the assumption  $\eta>3\beta/2$ and observe that the estimates we need,  still hold in the limiting case $\eta=3\beta/2$.  The analysis will require accurate estimates of the asymptotic behavior of the eigenvalues in the complex plane, with precise  estimates for the limiting case.
 
We will consider the following  Cauchy problem in the
square domain  $\Omega=(0,\pi)\times (0,\pi)$
\begin{equation}\label{id}
\begin{cases}
\displaystyle
u_{tt}(t,x,y) - \triangle u(t,x,y)+\beta\int_0^t\ e^{-\eta(t-s)} 
\triangle u(s,x,y)ds= 0\,,
\quad t\in (0,T),\,\, (x,y)\in \Omega,
\\
u(t,x,y)=0,\qquad   t\in (0,T), \,\, (x,y)\in \partial\Omega,
\\
u(0,x,y)=u_{0}(x,y),\quad
u_t(0,x,y)=u_{1}(x,y),\qquad  (x,y)\in \Omega.
\end{cases}
\end{equation}
In \cite{LS} we provided a detailed analysis of the cubic equation associated
with the integrodifferential equation. In particular, we gave the asymptotic behavior of the solutions of the cubic equation.
Using those results,  it is possible to write the solutions of the cubic equation in a different form with respect to that given in \cite{LS}, but more fitting for the goal of the present paper. Indeed, 
the following representation for the solution of problem \eqref{id} holds.
\begin{theorem}\label{th:repr}
For any  $(u_{0},u_{1})\in H^1_0(\Omega)\times 
L^2(\Omega)$ and $\eta\ge3\beta/2$ the weak solution of problem \eqref{id}  is given by
\begin{equation}\label{eqn:state}
u(t,x,y)=
\sum_{k_1,k_2=1}^\infty\Big(C_{k_1k_2}e^{ 
i\om_{k_1k_2}t}+\overline{C_{k_1k_2}}e^{-i\overline{\om_{k_1k_2}}t}
+R_{k_1k_2}e^{r_{k_1k_2}t}\Big)\sin(k_1 x)\sin(k_2 y),
\end{equation}
with
\begin{equation}\label{eq:rnI}
\begin{split}
\Re\omega_{k_1k_2}
&=\sqrt{k_1^2+k^2_2}\ \big(\Lambda^{-}_{k_1k_2}+\Lambda^{+}_{k_1k_2}\big),
\\
\Im\omega_{k_1k_2}&=\frac1{\sqrt3}\sqrt{k_1^2+k^2_2}\ \big(\Lambda^{-}_{k_1k_2}-\Lambda^{+}_{k_1k_2}\big)+\frac\eta3
\,,
\\
r_{k_1k_2}&=\frac2{\sqrt3}\sqrt{k_1^2+k^2_2}\ \big(\Lambda^{-}_{k_1k_2}-\Lambda^{+}_{k_1k_2}\big)-\frac\eta3\,,
\end{split}
\end{equation}
where 
\begin{equation}\label{fi}
\begin{split}
\Lambda^{-}_{k_1k_2}&={ 1\over 2}\root 3\of{\Phi_{k_1k_2}-\Psi_{k_1k_2}}\,,
\qquad
\Lambda^{+}_{k_1k_2}={ 1\over 2}\root 3\of{\Phi_{k_1k_2}+\Psi_{k_1k_2}}\,,
\\
\Phi_{k_1k_2}&=
\sqrt{1+
\Big(2\eta^2+{27\beta^2\over4}-9\eta\beta\Big)\frac1{k_1^2+k^2_2}
+{\eta^3(\eta-\beta)\over (k_1^2+k^2_2)^2}},
\\
\Psi_{k_1k_2}&=
{\eta^3\over 3{\sqrt{3(k_1^2+k_2^2)^3}}}
+\Big(\eta-{3\beta\over 2}\Big)\frac{\sqrt{3}}{\sqrt{ k_1^2+k_2^2}}\,.
\end{split}
\end{equation}
Moreover,
\begin{equation}\label{eqn:imkk}
r_{ k_1 k_2}\le -\Im\om_{k_1 k_2}\,,
\qquad 
\Im\om_{k_1 k_2}\le\frac\eta3
\qquad\qquad
\forall k_1,k_2\in\N,
\end{equation}and there exist $\mu>0$ such that
\begin{equation}\label{eqnRkk}
|R_{k_1 k_2}|\le\mu\frac{|C_{k_1 k_2}|}{\sqrt{ k_1^2+k_2^2}}
\qquad
\forall k_1,k_2\in\N
\,.
\end{equation}
\end{theorem}
\begin{remark}
In formula \eqref{eqn:state} the coefficients $C_{k_1k_2}$ and $R_{k_1k_2}$ are uniquely determined by the initial conditions $u_0$ and $u_1$. Since for our purposes  it is only significant  the relation \eqref{eqnRkk} between $R_{k_1k_2}$ and $C_{k_1k_2}$, we omit the explicit expression of $R_{k_1k_2}$.
\end{remark}
In virtue of Theorem \ref{th:repr} we are able to establish the following observability estimate on the subset $\Gamma=(0,\pi)\times\{0\}\cup\{0\}\times (0,\pi)$ of the boundary  of the
square domain.
 \begin{theorem}\label{eq:obsin}
Let  $\eta=3\beta/2$. 
If $u$ is the weak solution of problem \eqref{id} and 
$\Gamma=(0,\pi)\times\{0\}\cup\{0\}\times (0,\pi)$, then there exist $\beta_0>0$ and $T_0>0$ such that for all  $0<\beta<\beta_0$ and $T>T_0$ the inverse observability inequality 
\begin{equation}\label{eq:invobsin}
\int_0^T\int_{\Gamma}\Big|\frac{\partial 
u}{\partial\nu}\Big|^2d\Gamma dt
\ge
c_0\sum_{k_1,k_2=1}^\infty 
(k_1^2+k^2_2)|C_{k_1k_2}|^2(1+e^{-2\Im \om_{k_1k_2}T})
\,,
\end{equation}
holds true
for some positive  constant $c_0=c_0(T)$.
\end{theorem}

We will prove Theorem \ref{eq:obsin} in Section \ref{sub:proof} after  some preliminary results.

\section{Estimates of  the eigenvalues }
In this section we will study  the distribution of the eigenvalues in the complex plane.
Indeed, using the precise expressions of the eigenvalues provided by Theorem \ref{th:repr}, we will analyze the
behavior of partial gaps helpful to get the observability estimates.

To carry out our analysis, we need also the following known result, see \cite{KomLor(a)}. 
\begin{lemma}\label{le:2}
Fix an integer $N\ge 2$ and $N-1$ integers $k_1,\ldots,k_{N-1}\ge 1$.
If $k_N, k_N'$ are two positive integers satisfying
\begin{equation*}
\max\set{k_N, k_N'}\ge\max\set{k_1,\ldots,k_{N-1}},
\end{equation*}
then
\begin{equation*}
\abs{\sqrt{k_1^2+\cdots+k_{N-1}^ 2+k_N^2}-\sqrt{k_1^2+\cdots+k_{N-1}^ 2+(k_N')^2}}
\ge (\sqrt{N}-\sqrt{N-1})\abs{k_N-k_N'}.
\end{equation*}
In particular, if $N=2$ one has  $\sqrt{N}-\sqrt{N-1}=\sqrt{2}-1\approx 0.41$.

\end{lemma}

Using the notations introduced in Theorem \ref{th:repr} we will prove 
\begin{proposition}\label{p:repr}
If $\eta=3\beta/2$ and $\beta\in\big [0,\frac{2}{\sqrt{3}}\big]$, there exists $\gamma>0$ such that 
\begin{equation}\label{eq:k12}
\begin{split}
|\Re\omega_{k_1k_2}-\Re\omega_{k_1k'_2}|
&\ge\gamma |k_2-k'_2|
\qquad
\forall k_1\in\N,
\ \forall\max\set{k_2, k_2'}\ge k_1,
\\
|\Re\omega_{k_1k_2}-\Re\omega_{k'_1k_2}|
&\ge\gamma |k_1-k'_1|
\qquad
\forall k_2\in\N,
\ \forall\max\set{k_1, k_1'}\ge k_2,
\end{split}
\end{equation}
\begin{equation}\label{eq:reom}
\Re\om_{ k_1 k_2}\ge\gamma \sqrt{k_1^2+k^2_2}\,, \qquad\forall k_1,k_2\in\N.
\end{equation}
Moreover, the constant $\gamma=\gamma(\beta)$ in the previous inequalities can be taken equal to
\begin{equation}\label{eq:con-gamma}
\gamma=
\frac{\sqrt2-1}2\Bigg(\sqrt{1-\frac98 \beta^2
+\frac{27} {64}\beta^4}+\frac  {3\sqrt{3}}{16\sqrt2} \beta^3  \Bigg)^{1/3}
+\frac{\sqrt2-1}2\Bigg(\sqrt{1-\frac98 \beta^2
+\frac{27} {64}\beta^4}-\frac  {3\sqrt{3}}{16\sqrt2} \beta^3  \Bigg)^{1/3}
\,.
\end{equation}
In addition
\begin{equation}\label{eqn:imkk1}
0\le\Im\om_{k_1 k_2}\le\frac\beta2
\qquad
\forall k_1,k_2\in\N.
\end{equation}
\end{proposition}\label{pr:repr}
\begin{Proof}
First, by taking $\eta=\frac 32 \beta$ in formulas \eqref{eq:rnI} and \eqref{fi}, we obtain
\begin{equation}\label{eq:rnI1}
\Re\omega_{k_1k_2}
=\sqrt{k_1^2+k^2_2}\ \big(\Lambda^{+}_{k_1k_2}(\beta)+\Lambda^{-}_{k_1k_2}(\beta)\big)
\,.
\end{equation}
where
\begin{equation}\label{eq:phipsi}
\Phi_{k_1k_2}(\beta)=
\sqrt{1-\frac 9 4 \beta^2\frac1{k_1^2+k^2_2}
+\frac  {27} {16} \beta^4{1\over{(k_1^2+k^2_2)^2}}},
\qquad 
\Psi_{k_1k_2}(\beta)=
\frac  {3\sqrt{3}} 8 \beta^3 {1\over {\sqrt{(k_1^2+k_2^2)^3} } }\,,
\end{equation}
\begin{equation}\label{eq:Lambda+}
\begin{split}
\Lambda^{+}_{k_1k_2}(\beta)
&={ 1\over 2}\root 3\of{\Phi_{k_1k_2}(\beta)+\Psi_{k_1k_2}(\beta)}
\\
&=\frac12\Bigg(\sqrt{1-\frac 9 4 \beta^2\frac1{k_1^2+k^2_2}
+\frac  {27} {16} \beta^4{1\over{(k_1^2+k^2_2)^2}}}+\frac  {3\sqrt{3}} 8 \beta^3 {1\over {\sqrt{(k_1^2+k_2^2)^3} } } \Bigg)^{1/3}
\,,
\end{split}
\end{equation}
\begin{equation}\label{eq:Lambda-}
\begin{split}
\Lambda^{-}_{k_1k_2}(\beta)
&={ 1\over 2}\root 3\of{\Phi_{k_1k_2}(\beta)-\Psi_{k_1k_2}(\beta)}
\\
&=\frac12\Bigg(\sqrt{1-\frac 9 4 \beta^2\frac1{k_1^2+k^2_2}
+\frac  {27} {16} \beta^4{1\over{(k_1^2+k^2_2)^2}}}-\frac  {3\sqrt{3}} 8 \beta^3 {1\over {\sqrt{(k_1^2+k_2^2)^3} } } \Bigg)^{1/3}
\,.
\end{split}
\end{equation}
Fixed  $k_1\in\N$  and $k_2,k'_2\in\N$ with $k_2>k'_2$, thanks to \eqref{eq:rnI1} we have
 \begin{multline}\label{eq:Re1}
 \Re\omega_{k_1k_2}-\Re\omega_{k_1k'_2}
 =
\Big(\sqrt{k_1^2+k^2_2}- \sqrt{k_1^2+(k'_2)^2} \Big)\Big( \Lambda^{+}_{k_1k_2}(\beta)+\Lambda^{-}_{k_1k_2} (\beta)\Big)
\\
+\sqrt{k_1^2+(k'_2)^2}
\Big( \Lambda^{+}_{k_1k_2}(\beta)+\Lambda^{-}_{k_1k_2} (\beta)
- \Lambda^{+}_{k_1k'_2}(\beta)-\Lambda^{-}_{k_1k'_2} (\beta)\Big)
\,.
\end{multline}
We will show that the quantity $\Lambda^{+}_{k_1k_2}(\beta)+\Lambda^{-}_{k_1k_2} (\beta)$, regarded as function of $\frac1{\sqrt{k_1^2+k^2_2}}$, is decreasing for $\beta\in\big [0,\frac{2}{\sqrt{3}}\big]$. To do that, in view of \eqref{eq:Lambda+} and 
\eqref{eq:Lambda-}
we introduce the functions 
\begin{equation}\label{eq:ff}
F(x)=(f_+(x))^{1/3}+(f_-(x))^{1/3}\,,
\qquad
f_{\pm}(x)=\sqrt{1-x^2+\frac13x^4}\pm\frac{\sqrt3}9x^3,
\end{equation} 
since
\begin{equation}\label{eq:ff1}
\Lambda^{\pm}_{k_1k_2}(\beta)
=\frac12 f_{\pm}\Big(\frac{3\beta}{2\sqrt{k_1^2+k^2_2}}\Big)^{1/3}\,,
\qquad
\Lambda^{+}_{k_1k_2}(\beta)+\Lambda^{-}_{k_1k_2}(\beta)
=\frac12F\Big(\frac{3\beta}{2\sqrt{k_1^2+k^2_2}}\Big).
\end{equation}
We will prove that $F(x)$ is decreasing in $\big[0,\sqrt{\frac32}\big]$, that is $F'(x)\le0$ for any $x\in \big[0,\sqrt{\frac32}\big]$.
First, we note that
\begin{equation*}
F'(x)=\frac13(f_+(x))^{-2/3}f'_+(x)+\frac13(f_-(x))^{-2/3}f'_-(x)
\,.
\end{equation*}
Set $a(x)=\sqrt{1-x^2+\frac13x^4}$, in view of
\begin{equation*}
f_{\pm}'(x)=\frac{x}{3a(x)}\big(2x^2-3\pm a(x)\sqrt3x\big)
\qquad
x>0,
\end{equation*}
we can write
\begin{equation*}
F'(x)=\frac{x}{9a(x)}\bigg[(f_+(x))^{-2/3}(2x^2-3+ a(x)\sqrt3x\big)
+(f_-(x))^{-2/3}(2x^2-3- a(x)\sqrt3x\big)\bigg]
\,.
\end{equation*}
Therefore, $F'(x)\le0$ is equivalent to
\begin{equation*}
(f_+(x))^{-2/3}\big(a(x)\sqrt3x+2x^2-3\big)
\le(f_-(x))^{-2/3}\big(a(x)\sqrt3x-(2x^2-3)\big)
\,,
\end{equation*}
\begin{equation*}
(f_-(x))^{2}\big(a(x)\sqrt3x+2x^2-3\big)^3
\le(f_+(x))^{2}\big(a(x)\sqrt3x-(2x^2-3)\big)^3
\,,
\end{equation*}
and the last inequality is true for $f_+(x)>f_-(x)>0$ and $2x^2-3<0$, that is $x\in \big[0,\sqrt{\frac32}\big]$.
Therefore, it remains to be seen 
\begin{equation}\label{eq:f-pos}
f_-(x)>0
\qquad
\forall x\in \Big[0,\sqrt{\frac32}\Big]
\,.
\end{equation}
To this end, taking into account \eqref{eq:ff}, we note that  
 $\sqrt{1-x^2+\frac13x^4}>\frac{\sqrt3}9x^3$ if $x^6-9x^4+27x^2-27<0$.
Because of $x^6-9x^4+27x^2-27=(x^2-3)^3$ we have that $f_-(x)>0$ for $x<\sqrt3$.

%

Since for $\beta\in\big [0,\frac{2}{\sqrt{3}}\big]$ we have 
\begin{equation}\label{eq:f-pos1}
\frac{3\beta}{2\sqrt{k_1^2+k^2_2}}<\frac{3\beta}{2\sqrt{k_1^2+(k'_2)^2}}\le
\frac{3\beta}{2\sqrt{2}}\le\sqrt{\frac32}\,,
\end{equation}
thanks to \eqref{eq:ff1} we can deduce that
\begin{equation}\label{eq:decrease}
\Lambda^{+}_{k_1k_2}(\beta)+\Lambda^{-}_{k_1k_2} (\beta)
\ge\Lambda^{+}_{11}(\beta)+\Lambda^{-}_{11} (\beta)\,,
\qquad
\Lambda^{+}_{k_1k_2}(\beta)+\Lambda^{-}_{k_1k_2} (\beta)
\ge\Lambda^{+}_{k_1k'_2}(\beta)+\Lambda^{-}_{k_1k'_2} (\beta)\,,
\end{equation}
and hence, from \eqref{eq:Re1} it follows
\begin{equation}\label{eq:Re2}
 \Re\omega_{k_1k_2}-\Re\omega_{k_1k'_2}
 \ge
\big( \Lambda^{+}_{11}(\beta)+\Lambda^{-}_{11} (\beta)\big)
\Big(\sqrt{k_1^2+k^2_2}- \sqrt{k_1^2+(k'_2)^2} \Big)\,.
\end{equation}
Moreover, we also note that, thanks to $f_+(x)>-f_-(x)$, we have $F(x)>0$ for any $x\in\R$, whence we get 
$ \Lambda^{+}_{11}(\beta)+\Lambda^{-}_{11} (\beta)>0$.

\begin{figure}[H]
\begin{center}
 \includegraphics[]{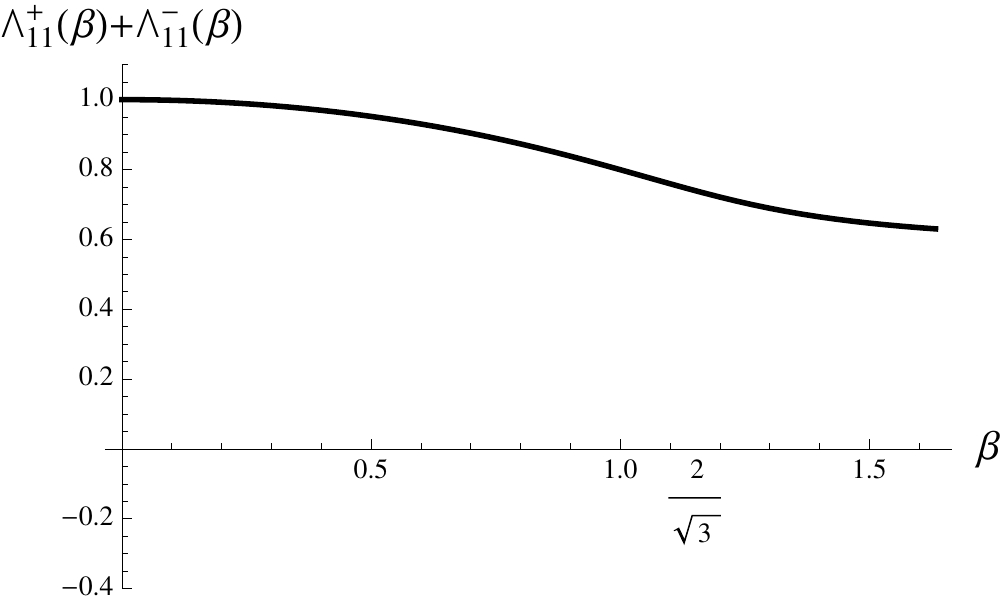}
\label{beta}
\end{center}
\end{figure}

Therefore from \eqref{eq:Re2}, using Lemma \ref{le:2},  we get
\begin{equation*}
\Re\omega_{k_1k_2}-\Re\omega_{k_1k'_2}
\ge (\sqrt2-1)\big( \Lambda^{+}_{11}(\beta)+\Lambda^{-}_{11} (\beta)\big)(k_2-k'_2),
\qquad
\forall k_1<k_2,
\end{equation*}
so, 
we have shown the first inequality in \eqref{eq:k12} with
$\gamma=(\sqrt2-1)\big( \Lambda^{+}_{11}(\beta)+\Lambda^{-}_{11} (\beta)\big)$, the same positive constant as in \eqref{eq:con-gamma}.
In a similar way one can prove the other inequality in \eqref{eq:k12}.

Finally,  in virtue of \eqref{eq:rnI1} and \eqref{eq:decrease}
we obtain
\begin{equation*}
\Re\om_{ k_1 k_2}
\ge\big( \Lambda^{+}_{11}(\beta)+\Lambda^{-}_{11} (\beta)\big)\sqrt{k_1^2+k^2_2}
\,, 
\qquad\forall k_1,k_2\in\N,
\end{equation*}
that is \eqref{eq:reom}.

Regarding the last statement \eqref{eqn:imkk1}, first we note that 
$\Im\om_{k_1 k_2}\le\frac\beta2$ follows from \eqref{eqn:imkk} since $\eta=3\beta/2$. 
To prove $\Im\om_{k_1 k_2}\ge0$, we have to take $\eta=\frac 32 \beta$ in \eqref{eq:rnI}  and show
\begin{equation*}
\Im\omega_{k_1k_2}
=\frac1{\sqrt3}\sqrt{k_1^2+k^2_2}\ \big(\Lambda^{-}_{k_1k_2}(\beta)-\Lambda^{+}_{k_1k_2}(\beta)\big)+\frac\beta2\ge0,
\end{equation*}
that is 
\begin{equation}\label{eq:posim}
\frac1{\sqrt3}\sqrt{k_1^2+k^2_2}\ \big(\Lambda^{+}_{k_1k_2}(\beta)-\Lambda^{-}_{k_1k_2}(\beta)\big)\le\frac\beta2,
\end{equation}
where $\Lambda^{+}_{k_1k_2}(\beta)$ and $\Lambda^{-}_{k_1k_2}(\beta)$ are given  by \eqref{eq:Lambda+} and \eqref{eq:Lambda-} respectively. We observe that
\begin{multline}\label{eq:pronot}
\Lambda^{+}_{k_1k_2}(\beta)-\Lambda^{-}_{k_1k_2}(\beta)
\\
=\frac{\Psi_{k_1k_2}(\beta)}{\root 3\of{(\Phi_{k_1k_2}(\beta)+\Psi_{k_1k_2}(\beta))^2}+\root 3\of{(\Phi_{k_1k_2}(\beta))^2-(\Psi_{k_1k_2}(\beta))^2}+\root 3\of{(\Phi_{k_1k_2}(\beta)-\Psi_{k_1k_2}(\beta))^2}}
\,.
\end{multline}
In virtue of \eqref{eq:ff1}, \eqref{eq:f-pos} and \eqref{eq:f-pos1} we have
\begin{equation*}
\root 3\of{\Phi_{k_1k_2}(\beta)-\Psi_{k_1k_2}(\beta)}
=2\Lambda^{-}_{k_1k_2}(\beta)
=f_{-}\Big(\frac{3\beta}{2\sqrt{k_1^2+k^2_2}}\Big)^{1/3}>0.
\end{equation*}
Therefore,  from \eqref{eq:pronot} we get
\begin{equation*}
\Lambda^{+}_{k_1k_2}(\beta)-\Lambda^{-}_{k_1k_2}(\beta)
\le\frac{\Psi_{k_1k_2}(\beta)}{\root 3\of{(\Phi_{k_1k_2}(\beta)+\Psi_{k_1k_2}(\beta))^2}}
\,,
\end{equation*}
whence,  taking also into account \eqref{eq:phipsi}, we have
\begin{equation*}
\frac1{\sqrt3}\sqrt{k_1^2+k^2_2}\ \big(\Lambda^{+}_{k_1k_2}(\beta)-\Lambda^{-}_{k_1k_2}(\beta)\big)
\le\frac1{\sqrt3}\sqrt{k_1^2+k^2_2}\ \root 3\of{\Psi_{k_1k_2}(\beta)}=\frac\beta2
\,,
\end{equation*}
that is  \eqref{eq:posim} holds true.
\end{Proof}

\begin{remark}
We observe that if  we pass to the limit in \eqref{eq:con-gamma} as $\beta\to0^+$, we obtain $\gamma=\sqrt2-1$, that is the value of the gap in the case of classical wave equations in a square domain, see Lemma {\rm\ref{le:2}}.
\end{remark}

\section{The observability estimate}
\subsection{The weight function}
To prove our results 
we need to introduce the  function
\begin{equation}\label{eq:k}
k(t):=\left \{\begin{array}{l}
\displaystyle\sin \frac{\pi t}{T}\,\qquad\qquad \mbox{if}\,\, t\in\ [0,T]\,,\\
\\
0\,\qquad\qquad\quad\  \ \ \  \mbox{otherwise}\,.
\end{array}\right .
\end{equation}
To begin with, we list some properties of $k$ in the following lemma.

\begin{lemma} \label{th:k}
Set
\begin{equation}\label{eqn:K}
K(u):=\frac{T\pi}{\pi^2-T^2u^2}\,,\qquad u\in \C\,,
\end{equation}
the following properties hold  for any $u\,,z\in \C$
\begin{equation*}
\overline{K(u)}=K(\overline{u})\,,
\end{equation*}
\begin{equation}\label{eq:sinek2bis}
\int_{-\infty}^{\infty} k(t)\Re(ze^{iu t})dt
= \Re\big(z(1+e^{iu T})K(u)\big)\,,
\end{equation}
\begin{equation}\label{eqn:sinek2}
\big|K(u)\big|=\big|K(\overline{u})\big|\,.
\end{equation}
For $\gamma>2\pi/T$, $j\in\N$ and $u\in\C$, $|u|\ge\gamma j$, we have
\begin{equation}\label{eq:sinek3}
\big|K(u)\big|\le
\frac{4\pi}{T\gamma^2(4j^2-1)}
\,.
\end{equation}
\end{lemma}
\begin{Proof}
We have to prove only the last statement. To this end,
we observe that
\begin{equation*}
\big|K(u)\big|=
\frac{\pi}{T\Big|u^2-\big(\frac{\pi}{T}\big)^2\Big|}
=\frac{4\pi}{T\gamma^2\Big|4\big(\frac{u}{\gamma}\big)^2-\big(\frac{2\pi}{T\gamma}\big)^2\Big|}
\,.
\end{equation*}
Since $| u|\ge\gamma j$ and $\frac{2\pi}{T\gamma}<1$,
we have
\begin{equation*}
\Big|4\Big(\frac{u}{\gamma}\Big)^2-\Big(\frac{2\pi}{T\gamma}\Big)^2\Big|
\ge
4\frac{|u|^2}{\gamma^2}-\Big(\frac{2\pi}{T\gamma}\Big)^2
\ge
4j^2-1\,,
\end{equation*}
and hence \eqref{eq:sinek3} follows.
\end{Proof}

\begin{theorem}\label{th:inverseest}
Assume that there exist $\gamma>0$ and 
$\tau\in\N$ such that
\begin{equation}\label{eqn:h1}
|\Re\om_n-\Re\om_m|\ge\gamma|n-m|
\qquad\forall n,m\in\N,\,\,\, \max\{n,m\}\ge\tau\,,
\end{equation}
and
\begin{equation}\label{eqn:h2}
\Re\om_n\ge\gamma n
\qquad\forall n\in\N\,,
\end{equation}
\begin{equation}\label{eq:rImom}
r_n\le -{\Im}\om_n
\qquad\forall\ n\in\N\,,
\end{equation}

\begin{equation}\label{eqn:rncn}
|R_n|\le\mu\frac{|C_n|}{n^\theta}
\qquad
\forall n\in\N
\qquad (\theta>1/2, \mu>0).
\end{equation}
Then
\begin{multline}\label{eq:inverseest}
\int_{0}^{T} \Big|\sum_{n=1}^\infty\Big(C_ne^{ 
i\om_nt}+\overline{C_n}e^{-i\overline{\om_n}t}+R_ne^{r_nt}\Big)\Big|^2\ 
dt
\\
\ge
2T\pi\sum_{n=\tau}^\infty\Big(\frac{1}{\pi^2+4T^2(\Im\om_n)^2}
-\frac{2 S}{T^2\gamma^2}\Big)|C_n|^2(1+e^{-2\Im \om_nT})
\\
-\frac{8\pi}{T\gamma^2 
}\Big(1+\frac{S}{2}\Big)\sum_{n=1}^\infty|C_n|^2(1+e^{-2\Im \om_n T})
\,,
\end{multline}
where $S=\mu\max\{\sum_{n=1}^{\infty}\frac{1}{n^{2\theta}},\frac{\pi^2}{6}\}$.
\end{theorem}
\begin{Proof}
Set
\begin{equation*}
F(t)=\sum_{n=1}^\infty\Big(C_ne^{ 
i\om_nt}+\overline{C_n}e^{-i\overline{\om_n}t}+R_ne^{r_nt}\Big),
\end{equation*}
we note that
\begin{multline*}
|F(t)|^2
=\Big|\sum_{n=\tau}^{\infty}\Big(C_ne^{ 
i\om_nt}+\overline{C_n}e^{-i\overline{\om_n}t}\Big)
+\sum_{n=1}^{\tau-1}\Big(C_ne^{ 
i\om_nt}+\overline{C_n}e^{-i\overline{\om_n}t}\Big)
+\sum_{n=1}^\infty R_ne^{r_nt}\Big|^2
\\
=\Big|\sum_{n=\tau}^{\infty}\Big(C_ne^{ 
i\om_nt}+\overline{C_n}e^{-i\overline{\om_n}t}\Big)\Big|^2
+2\sum_{n=\tau}^\infty\sum_{m=1}^{\tau-1}\Big(C_ne^{ 
i\om_nt}+\overline{C_n}e^{-i\overline{\om_n}t}\Big)
\Big(C_me^{ i\om_mt}+\overline{C_m}e^{-i\overline{\om_m}t}\Big)
\\
+2\sum_{n=\tau}^{\infty}\sum_{m=1}^\infty R_me^{r_mt}\Big(C_ne^{ 
i\om_nt}+\overline{C_n}e^{-i\overline{\om_n}t}\Big)
+\Big|\sum_{n=1}^{\tau-1}\Big(C_ne^{ 
i\om_nt}+\overline{C_n}e^{-i\overline{\om_n}t}\Big)
+\sum_{n=1}^\infty R_ne^{r_nt}\Big|^2
\,.
\end{multline*}
Let $k(t)$  be the function defined by (\ref{eq:k}). We have
\begin{multline*}
\int_{-\infty}^{\infty} k(t)|F(t)|^2\ dt
=
\int_{-\infty}^{\infty} k(t)
\Big|\sum_{n=\tau}^\infty\Big(C_ne^{ 
i\om_nt}+\overline{C_n}e^{-i\overline{\om_n}t}\Big)\Big|^2\ dt
\\
+2\int_{-\infty}^{\infty} k(t)
\sum_{n=\tau}^\infty\sum_{m=1}^{\tau-1}\Big(C_ne^{ 
i\om_nt}+\overline{C_n}e^{-i\overline{\om_n}t}\Big)
\Big(C_me^{ i\om_mt}+\overline{C_m}e^{-i\overline{\om_m}t}\Big)
\ dt
\\
+2\int_{-\infty}^{\infty} k(t)
\sum_{n=\tau}^\infty\sum_{m=1}^{\infty}R_me^{ r_mt}\Big(C_ne^{ 
i\om_nt}+\overline{C_n}e^{-i\overline{\om_n}t}\Big)
\ dt
\\
+\int_{-\infty}^{\infty} k(t)
\Big|\sum_{n=1}^{\tau-1}\Big(C_ne^{ 
i\om_nt}+\overline{C_n}e^{-i\overline{\om_n}t}\Big)
+\sum_{n=1}^\infty R_ne^{r_nt}\Big|^2
\ dt
\,.
\end{multline*}
Since $k(t)\ge0$ we can get rid of the last term on the right-hand 
side of the above formula, so we get
\begin{multline}\label{eqn:sum-s}
\int_{-\infty}^{\infty} k(t)|F(t)|^2\ dt
\ge
\int_{-\infty}^{\infty} k(t)
\Big|\sum_{n=\tau}^\infty\Big(C_ne^{ 
i\om_nt}+\overline{C_n}e^{-i\overline{\om_n}t}\Big)\Big|^2\ dt
\\
+2\int_{-\infty}^{\infty} k(t)
\sum_{n=\tau}^\infty\sum_{m=1}^{\tau-1}\Big(C_ne^{ 
i\om_nt}+\overline{C_n}e^{-i\overline{\om_n}t}\Big)
\Big(C_me^{ i\om_mt}+\overline{C_m}e^{-i\overline{\om_m}t}\Big)
\ dt
\\
+2\int_{-\infty}^{\infty} k(t)
\sum_{n=\tau}^\infty\sum_{m=1}^{\infty}R_me^{ r_mt}\Big(C_ne^{ 
i\om_nt}+\overline{C_n}e^{-i\overline{\om_n}t}\Big)
\ dt
\,.
\end{multline}
For $n,m\in\N$ we have
\begin{multline*}
\Big(C_ne^{ i\om_nt}+\overline{C_n}e^{-i\overline{\om_n}t}\Big)
\Big(C_me^{ i\om_mt}+\overline{C_m}e^{-i\overline{\om_m}t}\Big)
\\=
C_nC_me^{ i(\om_n+\om_m)t}+C_n\overline{C_m}e^{ i(\om_n-\overline{\om_m})t}
+\overline{C_n}C_me^{-i(\overline{\om_n}-\om_m)t}
+\overline{C_n}\overline{C_m}e^{-i(\overline{\om_n}+\overline{\om_m})t}
\\
=2\Re(C_n\overline{C_m}e^{ i(\om_n-\overline{\om_m})t}+C_nC_me^{ 
i(\om_n+\om_m)t})
\,,
\end{multline*}
so, by applying
  \eqref{eq:sinek2bis} we obtain
\begin{multline}\label{eq:lead}
\int_{-\infty}^{\infty} k(t)
\Big(C_ne^{ i\om_nt}+\overline{C_n}e^{-i\overline{\om_n}t}\Big)
\Big(C_me^{ i\om_mt}+\overline{C_m}e^{-i\overline{\om_m}t}\Big)\ dt
\\
=2\Re\big(C_n\overline{C_m}(1+e^{ 
i(\om_n-\overline{\om_m})T})K(\om_n-\overline{\om_m})
+C_nC_m(1+e^{ i(\om_n+\om_m)T})K(\om_n+\om_m)\big)
\,.
\end{multline}
Similarly,
by using again \eqref{eq:sinek2bis} we get
\begin{multline}\label{eq:misti}
\int_{-\infty}^{\infty} k(t)
e^{ r_mt}\Big(C_ne^{ i\om_nt}+\overline{C_n}e^{-i\overline{\om_n}t}\Big)
\ dt
=2\int_{-\infty}^{\infty} k(t)\Re\big(C_ne^{ i(\om_n-ir_m)t}\big)\ dt
\\
=2\Re\big(C_n(1+e^{ (i\om_n+r_m)T})K(\om_n-ir_m)\big)\,.
\end{multline}
By putting \eqref{eq:lead} and \eqref{eq:misti} into \eqref{eqn:sum-s}, we have
\begin{multline}\label{eqn:sum-s1}
\int_{-\infty}^{\infty} k(t)|F(t)|^2\ dt
\ge
\\
2\sum_{n,m=\tau}^\infty
\Re\Big(C_n\overline{C_m}(1+e^{ 
i(\om_n-\overline{\om_m})T})K(\om_n-\overline{\om_m})
+C_nC_m(1+e^{ i(\om_n+\om_m)T})K(\om_n+\om_m)\Big)
\\
+4\sum_{n=\tau}^\infty\sum_{m=1}^{\tau-1}
\Re\Big(C_n\overline{C_m}(1+e^{ 
i(\om_n-\overline{\om_m})T})K(\om_n-\overline{\om_m})
+C_nC_m(1+e^{ i(\om_n+\om_m)T})K(\om_n+\om_m)\Big)
\\
+4
\sum_{n=\tau}^\infty\sum_{m=1}^{\infty}R_m\Re\Big(C_n(1+e^{ 
(i\om_n+r_m)T})K(\om_n-ir_m)\Big)
\,.
\end{multline}
We may write the first sum on the right-hand side  as follows
\begin{multline*}
\sum_{n,m=\tau}^\infty
\Re\Big(C_n\overline{C_m}(1+e^{ 
i(\om_n-\overline{\om_m})T})K(\om_n-\overline{\om_m})\Big)
\\
=
\sum_{n=\tau}^\infty|C_n|^2 (1+e^{-2\Im \om_n T})K( \om_n-\overline{\om_n})
+
\sum_{\substack{n,m=\tau \\ n\not=m}}^\infty
\Re\Big(C_n\overline{C_m}(1+e^{ 
i(\om_n-\overline{\om_m})T})K(\om_n-\overline{\om_m})\Big)\,.
\end{multline*}
Plugging the above identity into (\ref{eqn:sum-s1})  we have
\begin{multline*}
\int_{-\infty}^{\infty} k(t)|F(t)|^2\ dt
\ge
2\sum_{n=\tau}^\infty|C_n|^2 (1+e^{-2\Im \om_n T})K( \om_n-\overline{\om_n})
\\
+2\sum_{\substack{n,m=\tau \\ n\not=m}}^\infty
\Re\Big(C_n\overline{C_m}(1+e^{ 
i(\om_n-\overline{\om_m})T})K(\om_n-\overline{\om_m})\Big)
+2\sum_{n,m=\tau}^\infty
\Re\Big(C_nC_m(1+e^{ i(\om_n+\om_m)T})K(\om_n+\om_m)\Big)
\\
+4\sum_{n=\tau}^\infty\sum_{m=1}^{\tau-1}
\Re\Big(C_n\overline{C_m}(1+e^{ 
i(\om_n-\overline{\om_m})T})K(\om_n-\overline{\om_m})
+C_nC_m(1+e^{ i(\om_n+\om_m)T})K(\om_n+\om_m)\Big)
\\
+4
\sum_{n=\tau}^\infty\sum_{m=1}^{\infty}R_m\Re\Big(C_n(1+e^{ 
(i\om_n+r_m)T})K(\om_n-ir_m)\Big)
\,.
\end{multline*}
By using the elementary estimate $\Re z\ge -|z|$, $z\in\C$, we obtain
\begin{multline}\label{eqn:sum-s2}
\int_{-\infty}^{\infty} k(t)|F(t)|^2\ dt
\ge 2\sum_{n=\tau}^\infty|C_n|^2(1+e^{-2\Im \om_n T})K( \om_n-\overline{\om_n})
\\
-2\sum_{\substack{n,m=\tau \\ n\not=m}}^\infty
|C_n| |C_m|(1+e^{-(\Im \om_n+\Im \om_m) T})\ |K( \om_n-\overline{\om_m})|
-2\sum_{n,m=\tau}^\infty|C_n| |C_m|(1+e^{-(\Im 
\om_n+\Im \om_m) T})\ |K( \om_n+\om_m)|
\\
-4\sum_{n=\tau}^\infty\sum_{m=1}^{\tau-1}
|C_n| |C_m|(1+e^{-(\Im \om_n+\Im \om_m) T})\Big( |K( 
\om_n-\overline{\om_m})|+|K(\om_n+\om_m)|\Big)
\\
-4\sum_{n=\tau}^\infty\sum_{m=1}^{\infty}|C_n|\ |R_m|(1+e^{(r_m-\Im 
\om_n) T}) |K( \om_n-ir_m)|
\,.
\end{multline}
By \eqref{eqn:sinek2} we have
$
|K( \om_n-\overline{\om_m})|=|K( \overline{\om_n}-\om_m)|\,,
$
and hence
\begin{equation}\label{eq:sum-s3}
\sum_{\substack{n,m=\tau \\ n\not=m}}^\infty
|C_n| |C_m|(1+e^{-(\Im \om_n+\Im \om_m) T})\ |K( \om_n-\overline{\om_m})|
\le
\sum_{n=\tau}^\infty|C_n|^2(1+e^{-2\Im \om_n T})\sum_{\substack{m=\tau \\ m\not=n}}^\infty\  |K( \om_n-\overline{\om_m})|
\,.
\end{equation}
Similarly
\begin{equation}
\sum_{n,m=\tau}^\infty
|C_n| |C_m|(1+e^{-(\Im \om_n+\Im \om_m) T})\ |K( \om_n+\om_m)|
\le
\sum_{n=\tau}^\infty|C_n|^2(1+e^{-2\Im \om_n T})\sum_{m=\tau}^\infty\ 
|K( \om_n+\om_m)|\,.
\end{equation}
Moreover,
\begin{multline}\label{eq:sum-s4}
\sum_{n=\tau}^\infty\sum_{m=1}^{\tau-1}
|C_n| |C_m|(1+e^{-(\Im \om_n+\Im \om_m) T})\Big( |K( 
\om_n-\overline{\om_m})|+|K(\om_n+\om_m)|\Big)
\\
\le
\frac12 \sum_{n=\tau}^\infty|C_n|^2(1+e^{-2\Im \om_n T})
\sum_{m=1}^{\tau-1}\Big( |K( \om_n-\overline{\om_m})|+|K(\om_n+\om_m)|\Big)
\\
+
\frac12\sum_{m=1}^{\tau-1} |C_m|^2(1+e^{-2\Im \om_m T})
\sum_{n=\tau}^\infty\Big( |K( \om_n-\overline{\om_m})|+|K(\om_n+\om_m)|\Big)
\,.
\end{multline}
Therefore, plugging formulas (\ref{eq:sum-s3})--(\ref{eq:sum-s4}) into 
(\ref{eqn:sum-s2}) we have
\begin{multline*}
\int_{-\infty}^{\infty} k(t)|F(t)|^2\ dt
\ge 2\sum_{n=\tau}^\infty|C_n|^2(1+e^{-2\Im \om_n T})K( \om_n-\overline{\om_n})
\\
-2\sum_{n=\tau}^\infty|C_n|^2(1+e^{-2\Im \om_n T})
\Big(\sum_{\substack{m=1 \\ m\not=n}}^\infty\  |K( \om_n-\overline{\om_m})|
+\sum_{m=1}^\infty\  |K( \om_n+\om_m)|\Big)
\\
-2\sum_{n=1}^{\tau-1} |C_n|^2(1+e^{-2\Im \om_n T})
\sum_{m=\tau}^\infty\Big( |K( \om_m-\overline{\om_n})|+|K(\om_m+\om_n)|\Big)
\\
-4\sum_{n=\tau}^\infty\sum_{m=1}^{\infty}|C_n|\ |R_m|(1+e^{(r_m-\Im 
\om_n) T}) |K( \om_n-ir_m)|
\,.
\end{multline*}
First, thanks to (\ref{eqn:K}) we note that
$$
K( \om_n-\overline{\om_n})=\frac{T\pi}{\pi^2+4T^2(\Im\om_n)^2}\,,
$$
so we can write
\begin{multline}\label{eqn:sum-s5}
\int_{-\infty}^{\infty} k(t)|F(t)|^2\ dt
\ge 2T\pi\sum_{n=\tau}^\infty|C_n|^2\frac{1+e^{-2\Im \om_n 
T}}{\pi^2+4T^2(\Im\om_n)^2}
\\
-2\sum_{n=\tau}^\infty|C_n|^2(1+e^{-2\Im \om_n T})
\Big(\sum_{\substack{m=1 \\ m\not=n}}^\infty\  |K( \om_n-\overline{\om_m})|
+\sum_{m=1}^\infty\  |K( \om_n+\om_m)|\Big)
\\
-2\sum_{n=1}^{\tau-1} |C_n|^2(1+e^{-2\Im \om_n T})
\sum_{m=\tau}^\infty\Big( |K( \om_m-\overline{\om_n})|+|K(\om_m+\om_n)|\Big)
\\
-4\sum_{n=\tau}^\infty\sum_{m=1}^{\infty}|C_n|\ |R_m|(1+e^{(r_m-\Im 
\om_n) T}) |K( \om_n-ir_m)|
\,.
\end{multline}
Now, for any $n\ge\tau$ we have to estimate the sum
$$
\sum_{\substack{m=1 \\ m\not=n}}^\infty\  |K( \om_n-\overline{\om_m})|
+\sum_{m=1}^\infty\  |K( \om_n+\om_m)|
\,.
$$
Since $\max\{n,m\}\ge\tau$, thanks to assumptions \eqref{eqn:h1} and 
\eqref{eqn:h2} we can apply
  (\ref{eq:sinek3}) to get
\begin{multline*}
\sum_{\substack{m=1 \\ m\not=n}}^\infty\  |K( \om_n-\overline{\om_m})|
+\sum_{m=1}^\infty\  |K( \om_n+\om_m)|
\\
\le
\frac{4\pi}{T\gamma^2 }
\Big(\sum_{\substack{m=1 \\ m\not=n}}^\infty\frac1{4(m-n)^2-1}
+\sum_{m=1}^\infty\frac1{4(m+n)^2-1}\Big)
\le
\frac{8\pi}{T\gamma^2 }
\sum_{j=1}^\infty\frac1{4j^2-1}
\,.
\end{multline*}
Since
$
\displaystyle
\sum_{j=1}^\infty\frac1{4j^2-1}=\frac12\,,
$
we have
\begin{equation*}
\sum_{\substack{m=1 \\ m\not=n}}^\infty\  |K( \om_n-\overline{\om_m})|
+\sum_{m=1}^\infty\  |K( \om_n+\om_m)|
\le
\frac{4\pi}{T\gamma^2 }
\,.
\end{equation*}
In view of the above estimate, we can write \eqref{eqn:sum-s5} in the 
following way
\begin{multline}\label{eqn:sum-s5bis}
\int_{-\infty}^{\infty} k(t)|F(t)|^2\ dt
\ge
2T\pi\sum_{n=\tau}^\infty|C_n|^2\frac{1+e^{-2\Im \om_n 
T}}{\pi^2+4T^2(\Im\om_n)^2}
-\frac{8\pi}{T\gamma^2 }\sum_{n=\tau}^\infty|C_n|^2(1+e^{-2\Im \om_n T})
\\
-2\sum_{n=1}^{\tau-1} |C_n|^2(1+e^{-2\Im \om_n T})
\sum_{m=\tau}^\infty\Big( |K( \om_m-\overline{\om_n})|+|K(\om_m+\om_n)|\Big)
\\
-4\sum_{n=\tau}^\infty\sum_{m=1}^{\infty}|C_n|\ |R_m|(1+e^{(r_m-\Im 
\om_n) T}) |K( \om_n-ir_m)|
\,.
\end{multline}
Concerning the third sum on the right-hand side of the previous 
estimate, as above we can show that
 \begin{equation*}
\sum_{m=\tau}^\infty\Big( |K( \om_m-\overline{\om_n})|+|K(\om_m+\om_n)|\Big)
\le
\frac{4\pi}{T\gamma^2 }
\qquad
\forall
n\le\tau-1
\,,
\end{equation*}
and hence 
\begin{equation*}
\sum_{n=1}^{\tau-1} |C_n|^2(1+e^{-2\Im \om_n T})
\sum_{m=\tau}^\infty\Big( |K( \om_m-\overline{\om_n})|+|K(\om_m+\om_n)|\Big)
\le
\frac{4\pi}{T\gamma^2 }
\sum_{n=1}^{\tau-1}|C_n|^2(1+e^{-2\Im \om_n T})
\,.
\end{equation*}
Therefore,  from \eqref{eqn:sum-s5bis} it follows
\begin{multline}\label{eq:sum-s5bis}
\int_{-\infty}^{\infty} k(t)|F(t)|^2\ dt
\ge
2T\pi\sum_{n=\tau}^\infty|C_n|^2\frac{1+e^{-2\Im \om_n 
T}}{\pi^2+4T^2(\Im\om_n)^2}
-\frac{8\pi}{T\gamma^2 }\sum_{n=1}^\infty|C_n|^2(1+e^{-2\Im \om_n T})
\\
-4\sum_{n=\tau}^\infty\sum_{m=1}^{\infty}|C_n|\ |R_m|(1+e^{(r_m-\Im 
\om_n) T}) |K( \om_n-ir_m)|\,.
\end{multline}
Now, we have to estimate the last sum on the right-hand side.
Thanks to (\ref{eqn:rncn}), we have
\begin{multline}\label{eqn:mistisin}
4\sum_{n=\tau }^\infty\sum_{m= 1}^{\infty}|C_n||R_m|\ |K( \om_n-ir_m)|
\le 4\mu
\sum_{n=\tau}^\infty\sum_{m= 1}^{\infty}|C_n|\frac{|C_m|}{m^\theta} \ 
|K( \om_n-ir_m)|\\
\le
2\mu\sum_{n=\tau}^\infty|C_n|^2\sum_{m= 1}^{\infty}\frac{|K( 
\om_n-ir_m)|}{m^{2\theta}}
+2\mu\sum_{m= 1}^{\infty}|C_m|^2\sum_{n=\tau}^\infty|K( \om_n-ir_m)|
\,.
\end{multline}
Since $\Re\om_n\ge\gamma n$, 
again by (\ref{eq:sinek3}) we have
\begin{eqnarray*}
|K( \om_n-ir_m)|
\le
\frac{4\pi}{T\gamma^2(4n^2-1)}\,.
\end{eqnarray*}
As a consequence,
we get
\begin{equation*}
2\mu\sum_{n=\tau}^\infty|C_n|^2\sum_{m= 1}^{\infty}\frac{|K( 
\om_n-ir_m)|}{m^{2\theta}}
\le
\frac{4\pi\mu}{T\gamma^2}\sum_{m=1}^{\infty}\frac{1}{m^{2\theta}}\sum_{n=\tau}^\infty\frac{|C_n|^2}{2n^2-1/2}
\le
\frac{4\pi\mu}{T\gamma^2}\sum_{m=1}^{\infty}\frac{1}{m^{2\theta}}\sum_{n=\tau}^\infty|C_n|^2
\,,
\end{equation*}
\begin{equation*}
2\mu\sum_{m= 1}^{\infty}|C_m|^2\sum_{n=\tau}^\infty|K( \om_n-ir_m)|
\le
\frac{4\pi\mu}{T\gamma^2}\sum_{n=1}^{\infty}\frac{1}{n^2}\sum_{m= 
1}^\infty|C_m|^2
\,.
\end{equation*}
Plugging the two previous estimates into \eqref{eqn:mistisin}, we obtain
\begin{equation*}
4\sum_{n=\tau}^\infty\sum_{m=1}^{\infty}|C_n||R_m|\ |K( \om_n-ir_m)|
\le
\frac{4\pi\mu}{T\gamma^2}
\sum_{n=1}^{\infty}\frac{1}{n^{2\theta}}\sum_{n=\tau}^\infty|C_n|^2
+\frac{4\pi\mu}{T\gamma^2}\sum_{n=1}^{\infty}\frac{1}{n^2}
\sum_{n= 1}^\infty|C_n|^2
\,.
\end{equation*}
In addition, keeping in mind also \eqref{eq:rImom}, in a similar way it follows
\begin{multline*}
4\sum_{n=\tau}^\infty\sum_{m= 1}^{\infty}|C_n||R_m|\ e^{(r_m-\Im 
\om_n) T} |K( \om_n-ir_m)|
\le
4\sum_{n=\tau}^\infty\sum_{m= 1}^{\infty}|R_m|e^{-\Im \om_mT}|C_n|\ 
e^{-\Im \om_nT} |K( \om_n-ir_m)|
\\
\le
\frac{4\pi\mu}{T\gamma^2}
\sum_{n=1}^{\infty}\frac{1}{n^{2\theta}}
\sum_{n= \tau}^\infty|C_n|^2e^{-2\Im \om_nT}
+\frac{4\pi\mu}{T\gamma^2}
\sum_{n=1}^{\infty}\frac{1}{n^2}
\sum_{n=1}^\infty|C_n|^2e^{-2\Im \om_nT}
\,,
\end{multline*}
and hence, by summing the previous two inequalities yields
\begin{multline*}
4\sum_{n=\tau}^\infty\sum_{m= 1}^{\infty}|C_n||R_m|\ (1+e^{(r_m-\Im 
\om_n) T}) |K( \om_n-ir_m)|
\\
\le
\frac{4\pi S}{T\gamma^2}
\sum_{n=\tau}^\infty|C_n|^2(1+e^{-2\Im \om_nT})
+\frac{4\pi S}{T\gamma^2}
\sum_{n= 1}^\infty|C_n|^2(1+e^{-2\Im \om_nT})
\,,
\end{multline*}
where 
$S=\mu\max\{\sum_{n=1}^{\infty}\frac{1}{n^{2\theta}},\frac{\pi^2}{6}\}$.
Finally, putting the previous formula in \eqref{eq:sum-s5bis}, we get the  estimate
\begin{multline*}
\int_{-\infty}^{\infty} k(t)|F(t)|^2\ dt
\ge
2T\pi\sum_{n=\tau}^\infty\Big(\frac{1}{\pi^2+4T^2(\Im\om_n)^2}
-\frac{2 S}{T^2\gamma^2}\Big)|C_n|^2(1+e^{-2\Im \om_nT})
\\
-\frac{8\pi}{T\gamma^2 
}\Big(1+\frac{S}{2}\Big)\sum_{n=1}^\infty|C_n|^2(1+e^{-2\Im \om_n T})
\,,
\end{multline*}
whence, in virtue of the definition of $k$, we obtain \eqref{eq:inverseest}
\end{Proof}
\subsection{Proof of Theorem \ref{eq:obsin}}\label{sub:proof}
Thanks to Theorem  \ref{th:repr} the weak solution of problem \eqref{id}
is given by
\begin{equation*}
u(t,x,y)=
\sum_{k_1,k_2=1}^\infty\Big(C_{k_1k_2}e^{ 
i\om_{k_1k_2}t}+\overline{C_{k_1k_2}}e^{-i\overline{\om_{k_1k_2}}t}
+R_{k_1k_2}e^{r_{k_1k_2}t}\Big)\sin(k_1 x)\sin(k_2 y).
\end{equation*}
Let $\Gamma_1=(0,\pi)\times\{0\}$. We have
\begin{multline}\label{eq:gamma1}
\int_0^T\int_{\Gamma_1}\Big|\frac{\partial u}{\partial\nu}\Big|^2d\Gamma dt
=\int_0^T\int_{0}^\pi|u_y(t,x,0)|^2dx dt
\\
=\int_0^T\int_{0}^\pi\Big|\sum_{k_1,k_2=1}^\infty
k_2\Big(C_{k_1k_2}e^{ i\om_{k_1k_2}t}+\overline{C_{k_1k_2}}e^{-i
\overline{\om_{k_1k_2}}t}
+R_{k_1k_2}e^{r_{k_1k_2}t}\Big)\sin(k_1 x)\Big|^2dx dt
\\
=\frac\pi2\sum_{k_1=1}^\infty\int_0^T\Big|\sum_{k_2=1}^\infty 
k_2\Big(C_{k_1k_2}e^{ 
i\om_{k_1k_2}t}+\overline{C_{k_1k_2}}e^{-i\overline{\om_{k_1k_2}}t}
+R_{k_1k_2}e^{r_{k_1k_2}t}\Big)\Big|^2 dt
\,.
\end{multline}
By Proposition \ref{p:repr} for any fixed $k_1$ the assumptions of Theorem  \ref{th:inverseest} are satisfied with $\tau=k_1$, so applying formula \eqref{eq:inverseest} we get 
\begin{multline*}
\int_{0}^{T} \Big|\sum_{k_2=1}^\infty k_2\Big(C_{k_1k_2}e^{ 
i\om_{k_1k_2}t}+\overline{C_{k_1k_2}}e^{-i\overline{\om_{k_1k_2}}t}
+R_{k_1k_2}e^{r_{k_1k_2}t}\Big)\Big|^2\ dt
\\
\ge
2T\pi\sum_{k_2=k_1}^\infty\Big(\frac{1}{\pi^2+4T^2(\Im\om_{k_1k_2})^2}
-\frac{2 S}{T^2\gamma^2}\Big)k_2^2|C_{k_1k_2}|^2(1+e^{-2\Im \om_{k_1k_2}T})
\\
-\frac{8\pi}{T\gamma^2 }\Big(1+\frac{S}{2}\Big)\sum_{k_2=1}^\infty 
k_2^2|C_{k_1k_2}|^2(1+e^{-2\Im \om_{k_1k_2} T}),
\end{multline*}
where
$S=\frac{\pi^2}{6}\mu$. The above formula can be also written in the following way
\begin{multline*}
\int_{0}^{T} \Big|\sum_{k_2=1}^\infty k_2\Big(C_{k_1k_2}e^{ 
i\om_{k_1k_2}t}+\overline{C_{k_1k_2}}e^{-i\overline{\om_{k_1k_2}}t}
+R_{k_1k_2}e^{r_{k_1k_2}t}\Big)\Big|^2\ dt
\\
\ge
2T\pi\sum_{k_2=k_1}^\infty\frac{1}{\pi^2+4T^2(\Im\om_{k_1k_2})^2}k_2^2|C_{k_1k_2}|^2(1+e^{-2\Im 
\om_{k_1k_2}T})
\\
-\frac{8\pi}{T\gamma^2 }\big(1+S\big)\sum_{k_2=k_1}^\infty 
k_2^2|C_{k_1k_2}|^2(1+e^{-2\Im \om_{k_1k_2}T})
-\frac{8\pi}{T\gamma^2 }\Big(1+\frac{S}{2}\Big)\sum_{k_2=1}^{k_1-1} 
k_2^2|C_{k_1k_2}|^2(1+e^{-2\Im \om_{k_1k_2} T})
\,.
\end{multline*}
Now, we note that in the previous estimate we may change $k^2_2$ into 
$|k|^2=k^2_1+k^2_2$. Indeed, if $k_2\ge k_1$ we have
$2k^2_2\ge|k|^2$, while for any $k_1\,, k_2$ we have
$k^2_2\le|k|^2$.
So, thanks also to \eqref{eqn:imkk1}, we obtain
\begin{multline*}
\int_{0}^{T} \Big|\sum_{k_2=1}^\infty k_2\Big(C_{k_1k_2}e^{ 
i\om_{k_1k_2}t}+\overline{C_{k_1k_2}}e^{-i\overline{\om_{k_1k_2}}t}
+R_{k_1k_2}e^{r_{k_1k_2}t}\Big)\Big|^2\ dt
\\
\ge
\frac{T\pi}{\pi^2+T^2\beta^2}\sum_{k_2=k_1}^\infty|k|^2|C_{k_1k_2}|^2(1+e^{-2\Im 
\om_{k_1k_2}T})
\\
-\frac{8\pi}{T\gamma^2 }\big(1+S\big)\sum_{k_2=k_1}^\infty 
|k|^2|C_{k_1k_2}|^2(1+e^{-2\Im \om_{k_1k_2}T})
-\frac{8\pi}{T\gamma^2 }\Big(1+\frac{S}{2}\Big)\sum_{k_2=1}^{k_1-1} 
|k|^2|C_{k_1k_2}|^2(1+e^{-2\Im \om_{k_1k_2} T})
\,.
\end{multline*}
Therefore, in view  of \eqref{eq:gamma1} it follows
\begin{multline}\label{eq:gamma11}
\int_0^T\int_{\Gamma_1}\Big|\frac{\partial u}{\partial\nu}\Big|^2d\Gamma dt
\\
\ge
\frac{T\pi^2}2\Big(\frac{1}{\pi^2+T^2\beta^2}
-\frac{8}{T^2\gamma^2 }\big(1+S\big)
\Big)\sum_{k_1=1}^\infty 
\sum_{k_2=k_1}^\infty|k|^2|C_{k_1k_2}|^2(1+e^{-2\Im \om_{k_1k_2}T})
\\
-\frac{4\pi^2}{T\gamma^2 }\Big(1+\frac{S}{2}\Big)
\sum_{k_1=1}^\infty\sum_{k_2=1}^{k_1-1} 
|k|^2|C_{k_1k_2}|^2(1+e^{-2\Im \om_{k_1k_2} T})
\,.
\end{multline}
In a similar way we can establish the following estimate for 
$\Gamma_2=\{0\}\times(0,\pi)$:
\begin{multline*}
\int_0^T\int_{\Gamma_2}\Big|\frac{\partial u}{\partial\nu}\Big|^2d\Gamma dt
\\
\ge
\frac{T\pi^2}2\Big(\frac{1}{\pi^2+T^2\beta^2}
-\frac{8}{T^2\gamma^2 }\big(1+S\big)
\Big)\sum_{k_2=1}^\infty 
\sum_{k_1=k_2}^\infty|k|^2|C_{k_1k_2}|^2(1+e^{-2\Im \om_{k_1k_2}T})
\\
-\frac{4\pi^2}{T\gamma^2 }\Big(1+\frac{S}{2}\Big)
\sum_{k_2=1}^\infty\sum_{k_1=1}^{k_2-1} 
|k|^2|C_{k_1k_2}|^2(1+e^{-2\Im \om_{k_1k_2} T})
\\
=
\frac{T\pi^2}2\Big(\frac{1}{\pi^2+T^2\beta^2}
-\frac{8}{T^2\gamma^2 }\big(1+S\big)
\Big)
\sum_{k_1=1}^\infty 
\sum_{k_2=1}^{k_1}|k|^2|C_{k_1k_2}|^2(1+e^{-2\Im \om_{k_1k_2}T})
\\
-\frac{4\pi^2}{T\gamma^2 }\Big(1+\frac{S}{2}\Big)
\sum_{k_1=1}^\infty\sum_{k_2=k_1+1}^\infty 
|k|^2|C_{k_1k_2}|^2(1+e^{-2\Im \om_{k_1k_2} T})
\,.
\end{multline*}
Thanks to the expression \eqref{eq:con-gamma} of $\gamma=\gamma(\beta)$, there exists $\beta_0>0$ such that $\gamma^2-8(1+S)\beta^2>0$ for any $\beta\in (0,\beta_0]$, and hence, for $T>2\pi\sqrt{\frac{2(1+S)}{\gamma^2-8(1+S)\beta^2}}$ we get
\begin{equation*}
\frac{1}{\pi^2+T^2\beta^2}
-\frac{8}{T^2\gamma^2 }\big(1+S\big)
=
\frac{T^2(\gamma^2-8(1+S)\beta^2)-8(1+S)\pi^2}{T^2\gamma^2(\pi^2+T^2\beta^2)}>0
\,.
\end{equation*}
Therefore, from the previous estimates we can deduce 
\begin{multline}\label{eq:gamma2}
\int_0^T\int_{\Gamma_2}\Big|\frac{\partial u}{\partial\nu}\Big|^2d\Gamma dt
\\
\ge
\frac{T\pi^2}2\Big(\frac{1}{\pi^2+T^2\beta^2}
-\frac{8}{T^2\gamma^2 }\big(1+S\big)
\Big)\sum_{k_1=1}^\infty 
\sum_{k_2=1}^{k_1-1}|k|^2|C_{k_1k_2}|^2(1+e^{-2\Im \om_{k_1k_2}T})
\\
-\frac{4\pi^2}{T\gamma^2 }\Big(1+\frac{S}{2}\Big)
\sum_{k_1=1}^\infty\sum_{k_2=k_1}^\infty 
|k|^2|C_{k_1k_2}|^2(1+e^{-2\Im \om_{k_1k_2} T})
\,.
\end{multline}
By summing \eqref{eq:gamma11} and \eqref{eq:gamma2} and taking into account that $\Gamma=\Gamma_1\cup\Gamma_2$ we get
\begin{equation*}
\int_0^T\int_{\Gamma}\Big|\frac{\partial 
u}{\partial\nu}\Big|^2d\Gamma dt
\ge
\frac{T\pi^2}2\Big(\frac{1}{\pi^2+T^2\beta^2}
-\frac{4}{T^2\gamma^2 }\big(4+3S\big)
\Big)\sum_{k_1,k_2=1}^\infty 
|k|^2|C_{k_1k_2}|^2(1+e^{-2\Im \om_{k_1k_2}T})
\,.
\end{equation*}
Finally, again in view of
\eqref{eq:con-gamma}, we can pick out $\beta_0>0$ sufficiently small such that $\gamma^2-4(4+3S)\beta^2>0$ for any $\beta\in (0,\beta_0]$, and hence 
$$
c_0:=\frac{T\pi^2}2\Big(\frac{1}{\pi^2+T^2\beta^2}-\frac{4}{T^2\gamma^2 }\big(4+3S\big)\Big)>0
\qquad
\forall T>2\pi\sqrt{\frac{4+3S}{\gamma^2-4(4+3S)\beta^2}}
\,,
$$
so
\eqref{eq:invobsin} holds true. $\square$

\end{document}